\documentclass[12pt]{article}
\usepackage[T2A]{fontenc}
\usepackage[cp1251]{inputenc}
\usepackage{amsmath}
\usepackage{amssymb}
\usepackage{amsthm}
\usepackage{graphicx}

\newtheorem{theorem}{Theorem}

\newtheorem{conjecture}{Conjecture}

\newtheorem{corollary}{Corollary}
\newtheorem{assumption}{Assumption}
\newtheorem{claim}{Claim}
\newtheorem{program}{Program}

\renewcommand{\div}{\mathop{\rm div}\nolimits}

\font\Bbb=msbm10 scaled 1200
\def\R{\hbox{\Bbb R}}

\def\O{\hbox{\Bbb O}}
\def\beq{\begin{equation}} \def\eeq{\end{equation}}

\begin{document}

523.46/.481
\begin{center}
          THE FRACTAL THEORY OF THE SATURN RING \\
                       M.I.Zelikin

%\footnotemark{} \footnotetext{Работа выполнена при финансовой
%поддержке РФФИ, грант № 11-01-00986-а, а также Программой
%Президиума РАН "Математическая теория управления".}
\end{center}

\bigskip

 \begin{abstract}
The true reason for partition of the Saturn ring as well as rings
of other planets into great many of sub-rings is found. This
reason is the theorem of Zelikin-Lokutsievskiy-Hildebrand about
fractal structure of solutions to generic piece-wise smooth
Hamiltonian systems. The instability of two-dimensional model of
rings with continues surface density of particles distribution is
proved both for Newtonian and for Boltzmann equations. We do not
claim that we have solved the problem of stability of Saturn ring.
We rather put questions and suggest some ideas and means for
researches.

 \end{abstract}

\section{Introduction}

Up to the middle of XX century the astronomy went hand in hand
with mathematics. So, it is no wonder that such seemingly purely
astronomical problems as figures of equilibrium for rotating,
gravitating liquids or the stability of Saturn Ring came to the
attention of the most part of outstanding mathematicians. It will
be suffice to mention Galilei, Huygens, Newton, Laplace,
Maclaurin, Maupertuis, Clairaut, d'Alembert, Legendre, Liouville,
Jacobi, Riemann, Poincar\'e, Chebyshev, Lyapunov, Cartan and
numerous others (see, for example, \cite{A}). It is pertinent to
recall the out of sight discussion between Poincar\'e and Lyapunov
regarding figures of equilibrium for rotating, gravitating
liquids. Poincar\'e, who obtained his results using not so
rigorous reasoning and often by simple analogy, wrote, ``It is
possible to make many objections, but the same rigor as in pure
analysis is not demanded in mechanics''. Lyapunov had entirely
different position: ``It is impermissible to use doubtful
reasonings when solving a problem in mechanics or physics (it is
the same) if it is formulated quite definitely as a mathematical
one. In that case, it becomes a problem of pure analysis and must
be treated as such''. Both great scientists had in a sense good
reasons. When it comes to real natural processes, one cannot be
aware that all external and interior effects were taken into
account. And so, any model is unusable as a formal proof that the
real process will behave in one or another way. As for figures of
equilibrium for rotating, gravitating liquids, both mathematicians
did not take into account its non-homogeneity (in contrast to
Clairaut), neglected of interior currents (in contrast to Riemann
and Dedekind), ignore electro-magnetic phenomena. No one, so far
as we know, tries to consider this last effect on figures of
equilibrium, though no doubt a close relations of electro-magnetic
phenomena with the gravitation must exist. It is suffice to
mention the electro-magnetic dynamo of Earth core. So, the refusal
of Poincar\'e to perform proofs with scrupulous attention is
psychologically justified. But on the other hand, the point of
view of Lyapunov, being more difficult, is more attractive for
mathematicians, because the essence of mathematics is the
abstract. If only one has an intuitive assurance that the main
cause is found, the purely formal proof is very much desirable. As
a matter of fact, the mathematical rigor is intimately associated
with the mathematical beauty.

\bigskip

Let us recall main facts about the history of explorations of the
Saturn Ring. When  Galilei directed his new made primitive
telescope to the heaven, he discovered, apart from the Jupiter
satellites, something strange related with Saturn. He could not
understand whether it was a triple joint ball or a ball with
handles. Hence, he encode his discovery by the famous anagram:
``Altissimum planetam tergeminum observavi'' which means: ``I have
observed the very remote Planet being triple''. A diversity of
conjectures had been put forward about this phenomenon. Galilei
himself considered Saturn as three different planets situated on
the same ray of sight or as the planet with two satellites and
sometimes called it "ears". When the Ring appears to vanish,
turning edgewise, he wondered: if Saturn as in the myth really had
swallowed his children. And only the great disciple of Galilei
--- Christian Huygens --- possessed courage to believe in his own
eyes and to behold something incredible: the ring of regular shape
hovering weightlessly in vacuity around the planet! Recently Rings
(though not so big as the Saturn Ring) were discovered around all
planets besides that of Earth-group: Jupiter, Uranus and Neptune
have Rings.

The first gap (the Cassini gap) of distribution of mass density in
the Ring was discovered in 1675. It has the width about 4800
kilometers. Later on a number of gaps with the width of hundreds
kilometers (Enke gap, Laplace gap, Huygens gap among others) were
detected. Situation was fundamentally changed after flights of
cosmic spacecrafts: Voyager 1 and 2, Pioneer, and NASA's Cassini
expedition. Photos that was made from these apparatuses allow to
calculate (it is difficult to believe!) until cent hundreds
concentric ringlets and the corresponding numbers of gaps and
divisions. The exact number is unknown because some ringlets stick
together and even intertwine to one another.

The theoretic investigations of the Saturn Ring began immediately
after its discovery. Researchers tried to find conditions that
enable the stability of the Ring. The stability as itself was of
no doubt being evident result of continues observations

But after discovery of the thin structure of the Ring, efforts of
scientists switch into explanations of the necessity of the
dividing the Ring on big number of ringlets. We will not speak on
groundless fantasies similar to the p-adic analysis, to the
quantization of orbits, to black matter and so on.

\bigskip

{\sl The point of view of this work is as follows: the exact
solution in the ideal case of the flat Ring that consists of
infinitesimal particles has to contain infinitely many ringlets.
This thesis is based on our theorem on generic fractal structure
of trajectories of piece-wise smooth Hamiltonian systems.
\cite{Z4}.}

\bigskip

There exists traditional distinction between three main
conjectures: the Ring is entirely rigid, or it is liquid or it
consists on many distinct particles (so called Cassini
conjecture).

Conjecture of rigid Ring was disproved by Laplace \cite{La}. He
proved that in case of constant density the rigid Ring would be
unstable. The arbitrary displacement of its center relative to the
gravity center of Saturn must lead to increasing of the anomaly,
and the Ring will collapse on Saturn. Hence such Ring has to be
inhomogeneous.

The remarkable work of Maxwell \cite{M}, that was honored by Adams
prize, anticipated all subsequent researches of the Saturn Ring on
one and a half centuries. At the beginning of his paper Maxwell
apologized to readers that he had explored the question, which is
not bounded with any practical benefit neither in navigation nor
in astronomy, being moved by scientific curiosity only. The method
of Maxwell consists in consideration of small oscillations
relative to the equilibrium state as solutions of variational
equations. Nowadays, as a result of progress in creation of modern
telescopes and especially due to flights of cosmic apparatuses to
Saturn, oscillating regimes became directly observable. One can
see spiral changes in density, bends of the disk, spokes, and so
on. In particular, Maxwell considered the question about the {\sl
nonhomogeneous} rigid Ring. Under the supposition that the
homogeneity is violated in a unique point where exists an
additional mass, he proved that it will be necessary for the
stability that this mass consists  as a minimum 4/5 of the total
mass of the Ring. But it means that actually it will be a
satellite with a comparatively small additional tail of asteroids,
that contradicts the observations.

Rado explored the case of continuous nonhomogeneous distribution
of masses. He proved that the necessary conditions for stability
is the following: the variation of the mass density has to be
changed from 0.04 to 2.7, that is no less than in hundred times
that again contradicts the observations. The Ring is extremely
relatively thin, and if some of its part was more thin that an
other, than (be Girn calculations) to sustain the satellites
attraction, its rigidity must be 1000 times as great as that of
the steel.

So the Ring cannot be rigid.

\bigskip

To explore the impact of the Ring itself on the behavior of its
particles Poincar\'e \cite{P} had found the main part of
asymptotic of the integral that describes the force of attraction
that acts on a point $P$ from a circle with the uniform linear
density. He introduced nice, new construction: the arithmetic-
geometrical mean. He proved that the main part of the asymptotic
of the attraction force for $P$ lying sufficiently near to the
circle is the arithmetic-geometrical mean of the nearest and the
farthest distant from $P$ to points of the circle. Poincar\'e used
this result in investigation of the model of the tore-like liquid
Ring with the ellipsoidal meridional section. Poincar\'e proved
its instability.

Maxwell \cite{M} \cite{M1} proved that the ultimate density of the
Ring which is compatible with the stability (for the {\sl liquid}
Ring as well as for the dusty ring) does not exceed 1/300 of the
density of Saturn. It is called the Maxwell limit. Poincar\'e
showed that if the {\sl liquid} Ring rotates with the speed of
satellite of the same orbit than its density must be greater than
1/16 of the planet density. It would be incompatible with the
Maxwell limit. So, the Ring cannot be liquid.

\bigskip

We are forced to accept Cassini conjecture: The Ring consists of
the set of distinct rigid asteroids. It counts in favor of this
conjecture that the Ring appears transparent; stars and the Saturn
surface are seen through it and the light passes the Ring without
refraction. The term itself "the Cassini conjecture" should be
considered as anachronism. Now it is not a conjecture but the
universally accepted model that is verified by repeated well
established facts of observational astronomy.

\bigskip
\bigskip

I remember a fantastic impression when I observes for the first time the Saturn Ring through a telescope. It was striking to see such a strange, excellent, and regular form among stars. At first, let us try to give a very rough, heuristic explanation of the phenomenon.

\bigskip

Saturn, having very big mass, attracts particles from the ambient space. Consider an auto-gravitating cloud of particles subject to mutual collisions in the field of gravity of Saturn. If we mentally switch out the gravitational force of Saturn and the process of colliding, then we obtain a flat disc $R < \bar R$ as a stable equilibrium \cite{A}. 
Now let us switch on the gravitational force of Saturn and the process of colliding. Let $\rho (R)$ be the density of particles, $l(R)$ be the length of free run.

\begin{claim}.

There exist positive constants $\rho^*$, $\varepsilon$ , and $R^*$ such that the density of particles generate an annulus

\beq
\rho (R) = \left \{\begin{array}{c}  \varepsilon \qquad \mbox{for } R<R^* \\ 
\rho^*- \varepsilon \qquad \mbox{for }\bar R <R < R^* \\
 \end{array}
\right.
\label{a}
\eeq

\end{claim}

Indeed. The probability of collision is growing together with the density. Hence, in the domain $R<R^*$, where the density is very low, the length of free run is big and collisions are scarce. The stochastic change of the angular
velocity is significant and the corresponding particles with the big probability fall into Saturn. So, the density of particles in the vicinity of Saturn has the tendency to fall.
On the contrary, in the domain $\bar R > R > R^*$, where the density of particles is high, collisions happen frequently. As a result the angular velocity attains its stationary value. The density changes negligible. In addition, the stock   
of particles is supplemented  from the ambient space due to gravitational force of Saturn.

So, in our rough approximation we obtain the step-function for the density of particles that generates the Ring of Saturn.

Let us remark that this argumentation relates not only to the Saturn Ring but to all big planet. Now it is known that Jupiter, Saturn, Uranus, and  Neptune have rings of its own.

\bigskip

{\sl In the present work it will be shown that the mass distribution in
the flat Ring consisted of the set of isolated gravitated
particles must be very irregular. It cannot be given by smooth
positive function. The statement is in accordant with the observed
decomposition of the Ring on numerous very small ringlets.}

\bigskip

\section{The model of flat Ring (Model $A$)}

\bigskip

Let us give some facts for convenience of readers. Particles of
the Ring consist in general of snow and ice with minor amount of
carbon and silicon. The size of particles varies from microns to
meters (very rare till kilometer). Diameter of the Ring is
approximately 480 000 km. Its thickness from 1 to 1.5 km. The mass
of Saturn is ${\cal M_S} \approx 6\cdot 10^{29} g$. The mass of
the Ring is ${\cal M_R} \approx 3\cdot 10^{22} g$.  The mass of
the Mimas (one of the nearest Moon of Saturn) is ${\cal M_M}
\approx 3\cdot 10^{19} g$.

\bigskip

The thickness of the Ring is insignificant in comparison with its
size. Hence, the approximation of the flat Ring seems reasonable.
Let us give the following reasoning as an additional argument.
Take the carrying plane of the Ring as the plane $Ox_1x_2$ and its
normal in the center of the Ring take as the axis $Ox_3$.
Particles, that for some reason appears out of the plane
$Ox_1x_2$, are subject of action of the additional attractive
force from the Ring. If the particle appears sufficiently far from
the edge of the Ring end sufficiently near to its plane than the
vertical component of the attraction force induced by the Ring is
practically constant. It is approximately equal to the force of
attraction of infinite plane $Ox_1x_2$. Hence, particles leaving
the plane $Ox_1x_2$ in addition to the movement under the action
of the central force execute harmonic oscillations directed by the
axis $Ox_3$. After averaging over these frequencies, the
dependence on $x_3$ disappears and it remains only along the plane
$Ox_1x_2$. This heuristic argumentation shows that two-dimensional
model has to be adequate to the situation.

So, at first we consider a collision-free two-dimensional model
(call it Model $A$) for particles lying in the round annulus $Q$
with the centrum at the point $O$. The pare of cartesian
coordinates $(x_1,x_2)$ will be denoted in what follows by a
single capital letter $X$. In the plane $Ox_1x_2$ we introduce
polar coordinates $R,\phi$ in addition to the cartesian
coordinates. Let the internal board of the Ring be $R=R_1$; the
external one be $R=R_2$.

\begin{assumption}

Suppose that a closed, two-dimensional annulus $Q$ is filled by
particles with the surface density $f^0(R)$ that is symmetric
relative to the group of rotations $\O_2$. Let the support of the
function $f_0$ be $Q$. Suppose that the density is continues and
positive at points of $Q$. It follows from the compactness of $Q$
that

 \beq
 f^0(X) > \alpha > 0
 \label{1}
 \eeq
\noindent for $X \in Q$.
\end{assumption}

Let us find the main member of the asymptotic of attraction forces
of the Ring which act on points lying near its boarder.

\bigskip

Consider a domain $\Omega \subset \R^2$ with smooth boundary and
with the density $f(X)$ of mass distribution. The resultant force
of attraction from the mass which contained in $\Omega$ acting on
a sampling point $X$ will be denoted by $F(X)$. The constant of
gravitation we take as 1.

\begin{claim}

Let the function $f(X)$ in a neighborhood of a point $X_0 \in
\partial \Omega$ meets conditions of the assumption 1. Than the
main part of the asymptotic of attraction forces $F(X)$ as $X \to
X_0 \in \partial\Omega$ is equivalent to $- C\ln |X-X_0|$.
\end{claim}

Proof.

Let the point $X_0$ be $O$, and consider the domain
$\Omega_\epsilon$ with the mass density $f(x)$, defined by
inequalities $\Omega_\epsilon = \{-A\le x_1 \le A, \; \epsilon \le
x_2 \le B\}$. The sign of $\epsilon$ is arbitrary. Let us estimate
the main part of the attraction force acting on a point $O$ from
the domain $\Omega_\epsilon$. In view of (\ref{1}), one has

 \beq
\int_\epsilon^B dx_2\int_{-A}^A \frac
{x_2f(x_1,x_2)}{(x_1^2+x_2^2)^{3/2}}dx
> 2\alpha \int_0^A\left \{\frac {1}{(x_1^2 + \epsilon^2)^{1/2}} -
\frac {1}{(x_1^2 + B^2)^{1/2}} \right \}dx_1 \approx
 C\ln \frac {1}{|\epsilon|}.
 \label{2}
 \eeq

$\Box$

Let us denote by ${\cal F}(R)$ the lamp gravitational force acting
on a particle from Saturn and its Ring. Equations of movement in
the model $A$ are

\beq
 \left \{
 \begin{array}{c}
\dot X = V \\
\dot V = {\cal F}(R) \\
 \end{array}
 \right .
\label{3}
 \eeq

\begin{theorem}

Let the assumption 1 be fulfilled. Than the state that is defined
by the function $f^0(R)$ is unstable.
\end{theorem}

Proof.

The additional force of gravitation, that was induced by the Ring
itself, acting at points lying far from the edge of Ring, is small
comparatively to that acting from Saturn. The reason is that the
density of the Ring does not exceed 1/300 of the density of
Saturn. Besides, for points $P$ that is placed inside of the Ring,
the force from the inner (relative to $P$) part of the Ring is
almost balanced by the outer part. Meantime, the formula (\ref{2})
follows that the Ring attraction infinitely grows as $P$ tends to
the boarder of the Ring.

Due to the $\O_2$-symmetry, it is necessary for $f^0(R)$ to be
steady-state that the function {\cal F}(R) does not change along
trajectories. Hence, particles have to move along circles. But for
particles with $X_0 \in \partial Q$ this movement is demanded the
infinite speed that is physically non-realizable. So, in the real
movement, particles would move inside the Ring (both for points
lying at $R=R_1$ and that lying at $R=R_2$). Hence if the Ring
appears in the state that meets the assumption 1, than particles
of the outer part must tend to approach Saturn, while that of the
inner part to withdraw from Saturn. As a result the Ring as a
whole should begin to contract. This process will proceed till the
approximation of the flat Ring becomes inadequate. Therefore the
assumption 1 leads to instability.

$\Box$

\smallskip

Though it is not directly connected with the proof of the theorem,
it will be interesting to follow the evolution of the system with
the distribution $f^0(t_0,R)$. Taking into account that the total
force from Saturn and from the Ring is directed by the radius, one
can draw the graph of the force ${\cal F}(R)$, that is shown
schematically on the Fig 1 (the exact form of the graph
depends on the distribution function $f$).

\begin{figure}
	\centering
	\includegraphics{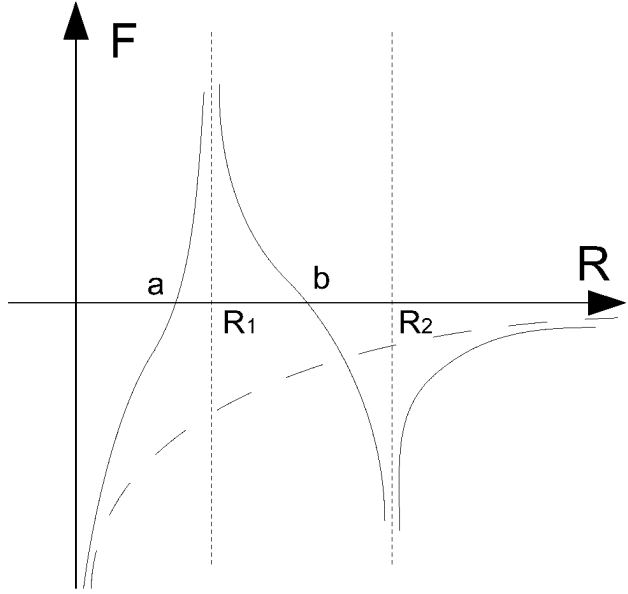}
	\caption{Attracting forces of the Ring and Saturn}
	\label{fig1}
\end{figure}

In the simplest case when the distribution function is constant on
the interval $[R_1,R_2]$, the function ${\cal F}(R)$ monotonically
increases from $-\infty$ to $+\infty$ on the interval $(0,R_1)$;
monotonically decreases from $+\infty$ to $-\infty$ on the
interval $(R_1,R_2)$ and monotonically increases from $-\infty$ to
zero on the interval $(R_2,+\infty)$.

In this case there exist exactly two points of intersection of the
graph of the function ${\cal F}(R)$ with the absciss axis: $R=a$
and $R=b$. These points correspond to libration lines. The name
libration means that these lines consist of fixed points of the
system (\ref{3}). It follows from consideration of signs of the
function ${\cal F}(R)$ that the circle $R=b$ consists of stable
equilibrium while the circle $R=a$ consists of unstable ones. By
using the theorem about continues dependance on initial data for
the system (\ref{3}), one can prove that the density $f$ remains
to meet conditions of assumption 1 in the process of growing of
$t$. Indeed, the circular symmetry will be conserved and the
strictly positivity on the interval $[R_1(t),R_2(t)]$  follows
from the matching of the density in corresponding under the
transference points and from the theorem about the integral
invariant for area. Hence, the function $R_1(t)$ monotonically
grows and the function $R_2(t)$ monotonically decreases. The Ring
will contract to the direction of the circle $R=b$. But we have to
take into account that the value $b$ itself will change in this
process.

One can put in doubt the existence of nonzero limit of the density
at the border of the Ring. But the observations show, that the
Ring and each of his ringlets have the precise, accurate
boundaries, instead of fuzzy boarder.

\bigskip

\begin{corollary}
The theorem 1 remains valid for any sufficiently wide sub-ring
which is divided from other parts of the Ring by two divisions
that is free from gravitating particles.
\end{corollary}

In the model $A$ we did not take into account effect of collision
of particles which a priory can substantially change the
situation. We have to consider the kinetic Boltzmann equation.
Consider the corresponding model $B$ when effect of collisions is
commensurable with the dynamical ones.

We have to modify a little the Boltzmann equation to apply it to
the two-dimensional situation, to the rotation of particles in the
Ring, to take into account the distribution of mass of particles,
and lows of its scattering.

\bigskip

\section{Kinetic Boltzmann equation (model $B$)}

\bigskip

To explore the stability of the Ring with regard to collisions we
have to apply methods of statistical physics and ergodic theory
\cite{G}, \cite{D}. Consider the integro-differential kinetic
Boltzmann equation \cite{B}, \cite{S}. In 1932 Carleman \cite{Ca}
proved the first global existence theorem for Boltzmann equation
in the space-homogeneous case. In 1963 Grad \cite{Gr} proved the
general local existence theorem in the neighborhood of the Maxwell
distribution. Later on this theorem was repeatedly generalized
(see, for example, the bibliography in \cite{U}).

The Boltzmann equation in connection with the Saturn Ring was
considered in works of Friedman and Gorkavyj \cite{GF}. Their
predecessors sought reasons of stability of the Ring because it
does not destruct almost 4 century. But Friedman and Gorkavyj
began to seek reasons of instability, namely, to seek oscillation
processes making the Ring to stratify in big number of thin
ringlets while preserving it as a whole. They assume that
stratification of the Ring is summoned by non-equilibrium
wave-like steady processes which take place in the Ring.

The kinetic Boltzmann equation plays the leading part in the
thermodynamics, the hydrodynamics, and the quantum mechanics. In
particular, from this equation follow Navier-Stokes equation,
Euler equation, Vlasov equation for plasma... (see, for example,
\cite{D}, \cite{Po}, \cite{U}).

To use Boltzmann equation one needs that the time of free running
essentially exceeds the time of colliding. This condition is
fulfilled in our situation. Indeed, mutual relative speeds between
particles of the Ring are estimated by values of several cm/cek.
The Ring is transparent, that means that sufficiently amount of
rays of the light do not meet particles of matter on his way
through the Ring, i.e. on a distant about 1-1.5 km. Besides,
estimations of the density and of mean size of particles show that
the time of free running are estimated by several hours. So, the
main conjecture about the time of collisions which is the base of
derivation of the Boltzmann equation may be thought of as
fulfilled.

\bigskip

Let us denote by $f_\mu (t,X,V)$ the scalar density of
distribution of particles with the mass $\mu$ in a moment $t$ at a
point $X \in Q$ that have the velocity $V$. Scalar product of
functions from the velocity $V$ will be defined by the metric
$L_2$ and will be designed by angular brackets. For instance, the
mean value of the function $\phi (V)$ for the density $f(V)$
equals

$$
<\phi (V),f(V)> = \int \phi (V)f(V)dV.
$$

Let us introduce the main macroscopic values, defined by a
distribution $f$. We use here the standard terminology of the
kinetic theory of gases. Note that the term "microscopic value",
being applied to particles with size around a meter (and even
reach amount a kilometre), seems at any rate strange. But one
cannot forget about astronomical size of the Ring itself.

The macroscopic density of particles equals

\beq
 \rho (t,X) =<1,f> = \int f(t,X,V)dV.
 \label{4}
 \eeq

The mathematical expectation of velocity per unit of mass is
defined by the vector $U$ with coordinates

 \beq
 u_j(t,X) = \frac {1}{\rho (t,X)} \int v_jf(t,X,V)dV.
 \label{5}
 \eeq

Denote by $W(t,X) = V - U(t,X)$ the chaotic velocity of particle
(the deviation of the real velocity from its mathematical
expectation).

Kinetic temperature per unit of mass equals

 \beq
 T(t,X) = \frac {1}{2\rho (t,X)\kappa_B} \int W^2(t,X)f(t,X,V)dV,
 \label{6}
 \eeq
where $\kappa_B$ is the Boltzmann constant. The term kinetic
temperature does not imply neither the temperature of ambient
space nor that of inside of particles. The kinetic temperature is
defined (as in kinetic theory of gases) as dispersion of a random
field of velocities of particles contained in the Ring.

The kinetic pressure is defined by the formula

 $$
  p = \kappa_B \rho T,
$$

The last relation is called by the equation of states in kinetic
theory of gases.

Unlike of the classical deduction of the Boltzmann equation we
have to take into account that particles can have all a priory
admissible masses. Hence one has to consider a stochastic
distribution $s(m)$ of masses. As a pattern we take the Boltzmann
arguments \cite{B} who considers the mixture of two gases: one
consists of molecules of the mass $m_1$ and other of the mass
$m_2$. Since the Boltzmann reasoning relates to the single act of
collision, it can be applied to the mixture of gases with any
given mass distribution, by introducing continuum distribution
functions to suit the mass of particles.

Given the equations of particles dynamics (\ref{3}) and the rules
for the changes its velocities after collisions one can write out
the equation for change of the density.

So, the kinetic Boltzmann equation is the balance equation for
change of the particles density of various masses subject to
forces ${\cal F}(X)$ and subject to changes of the movement due to
collisions.

\beq
 \frac{\partial f_\mu(t,X,V)}{\partial t} +
 \sum_{i=1}^2 \left ( \frac{\partial f_\mu(t,X,V)}{\partial x_i}v_i
 \right ) +  \sum_{i=1}^2 \left (
 \frac{\partial f_\mu(t,X,V)}{\partial v_i}{\cal F}_i(X)
 \right ) = J[f],
 \label{7}
 \eeq
where $J[f]$, that is called the collisions integral, is a
functional from the distribution function $f$ having the form

 \beq
 J[f] =  \frac {1}{\lambda}\int_0^\infty s(m)dm\int K(\sigma)
 [\tilde f_\mu \tilde f_m - f_\mu f_m]d\omega_\mu d\sigma =0.
  \label{8}
 \eeq

Here $d\omega_\mu$ is the volume element at the point, where the
particle of the mass $\mu$ lays; the kernel of the integral
$K(\sigma)$ describes the result of collision of two particles;
the value $\lambda$ is the time of free run. The sign tilde over
letters (say, $\tilde f$) shows on values of distribution function
at points that relate to the velocities after collision.

As it is known, the unique stationary (does not depend on time)
stable solution to the Boltzmann equation in $d$-dimensional space
is obtained only when the particles velocity have the Maxwell
distribution.

\beq M_{(\rho ,U,T)}(V) = \frac{\rho}{(2\pi \kappa_B T)^{d/2}}\exp
\left (- \frac{|V-U|^2}{2\kappa_B T}  \right ).
 \label{9}
 \eeq

In the two-dimensional case the exponent $d/2$ that stands into
denominator equals 1. The particles density distribution relative
to the masses is a result of subdivisions and sticking together of
particles in the colliding process. Parameters characterizes
collisions depend on relative velocities of colliding particles.
So, it is natural to suppose that the distribution of mass
particles in the process of relaxation of equilibrium with the
passage of time get adjusted to the velocity distribution and
become Maxwell-like.

The case is reduced to the addition into the Boltzmann equation
one more variable related to mass and to the introduction the set
of distribution functions which relate to each value of masses.
One has only to add under the integral sign the same factor as
that for the Maxwell distribution relative to velocities. Hence,
integrals of collisions remain in practice the same Gaussian
integrals as for the one-particle equation.

Collision integrals for the Maxwell distribution $M$ vanish and
the system (\ref{7}) that consists of continuum equations
decomposes into independent of one another scalar equations.

Methods of exploration of non-stationary solutions to the
Boltzmann equation were developed in \cite{C}, \cite{Ca},
\cite{Gr}, \cite{Ki}, \cite{U}, \cite{L} a.o.. The basic method is
to expand the unknown distribution function in a series (as in
perturbation theory) by moments of some fixed balanced
distribution $f$, as a rule, of Maxwell distribution $M$. The
benefit of such approach is that the calculation of moments needs
only the finite number of the expansion coefficients (its number
is equal to the degree of the moment).

\begin{assumption}

Suppose that forces ${\cal F}$ do not depend on time, and depend
only on the distance $R$ from the origin, and are defined by
potential of the field of Saturn and the Ring. Suppose that the
space distribution of particles is given by a smooth, positive
density. Suppose that the speed distribution $f^0_\mu (\cdot
,\cdot ,V)$ and the corresponded mass distribution both are
locally Maxwell-like.
\end{assumption}

The term "locally Maxwell" in this situation means that parameters
of Maxwell distribution (the mean value and the dispersion) depend
only on the distance from the origin.

It was shown above that if one does not take into account
collisions than the corresponding movement is unstable. Now we are
in the frame of  the thermodynamics and the kinetic theory when
the movement is described by the Boltzmann equation (\ref{7}).
This model can be regarded as an analog of the shear
two-dimensional Couette flow \cite{Co}, \cite{L}, \cite{R},
\cite{Sh}, when parallel walls that carry away the fluid have
constant but differing velocities. In that case the fluid is
stratified on fibers with constant velocities that vary linearly
from one fiber to another. The difference consists in the
following: in our case the sheared in the radius rotation is
caused by the external forces of gravitation that depend on the
radius. To describe this situation it is used the term
differential rotation.

So, we have two-dimensional model (call it as model $B$) of the

Ring $Q$ with the differential (depending on $R$) rotation being
induced by the gravitational forces of Saturn and the Ring.

Some theoretic arguments and experimental photometric data show
that ice particles of the Ring are covered by a thin loose
sediment that is made of snow. This layer soften strokes of
particles that take up a part of kinetic energy. This process
leads to equalization chaotic velocities. That means decreasing of
the dispersion i.e. the lowering of the kinetic temperature. But
there exists an opposite tendency that prevents cooling the Ring:
it is heating due to sheared viscosity, i.e. from the friction
between fibers with different speeds of rotation (as in the
Couetta flow). At the state of equilibrium shear heating must
compensate the dissipation of the energy in collisions.

\begin{theorem}

The stationary steady state solution to the system (\ref{7}) of
the Boltzmann equation for the model $B$ with free boundary
conditions cannot be given by distribution functions
$f_\mu^0(X,V)$ which meet assumptions 1 and 2.
\end{theorem}

Let us recall that these assumptions mean that functions
$f_\mu^0(X,V)$ are $C^1$-smooth and bounded from below by a
positive constant in the domain $Q$. Secondly they  have local
Maxwell distribution with parameter depending only on $R$, and the
same is true for the mass distribution.

Proof.

We do not need in the complete using of methods of Chapman and
Enskog \cite{C} and to calculate the exact steady state
distribution. It will be suffice to show that such distribution
cannot meet conditions of the assumptions 1 and 2.

To describe the transport equation for macroscopic values
 $\rho (X,t), U(X,t), T(X,t)$ one has to multiple the Boltzmann
 equation (\ref{7}) by $1, V, T$ correspondingly and to integrate
 over $V$.

In the first case the result of this operation for the density of
particles gives

$$
\frac{\partial \rho (t,X)}{\partial t} +
 \int\left ( \sum_{i=1}^2 \frac{\partial f_\mu^0(X,V)}{\partial
 x_i}v_i  \right )dV +  \int \sum_{i=1}^2 \left (
 \frac{\partial f_\mu^0(X,V)}{\partial v_i}{\cal F}_i(X)
 \right ) dV = 0.
$$

The right  hand site is zero since $f^0_\mu(X,V)$ is the Maxwell
distribution. By the same reason, the last integral in the left
hand site equals zero. Indeed, after carrying out of the sign of
the integral the factor ${\cal F}_i(X)$ which does not depend on
$V$, it remains the integral from the full derivative of Gaussian
distribution that equals zero.

The integration by parts of the remaining integral gives equation
of continuity

\beq
 \frac{\partial \rho (t,X)}{\partial t} +
 \div (\rho (t,X)U(t,X))=0.
 \label{10}
 \eeq

In the second case one obtains equation for mathematical
expectation of velocity (\ref{5})

$$
\frac{\partial (\rho (t,X) u_j(t,X))}{\partial t} + \int
\sum_{i=1}^2 v_jv_i\frac {\partial f_\mu^0(X,V)}{\partial x_i}dV +
\int \sum_{i=1}^2 v_j{\cal F}_i(X) \frac {\partial
f_\mu^0(X,V)}{\partial v_i}dV = 0.
$$

Again, after carrying out of the sign of the integral factors
which do not depend on $V$ and after integrating by parts, one
obtains

$$
\frac{\partial (\rho (t,X)u_j(t,X))}{\partial t} + \sum_{i=1}^2
\frac {\partial}{\partial x_i} \left \langle v_iv_j, f_\mu^0(X,V)
 \right \rangle - \rho (t,X) {\cal F}_j(X) = 0.
$$

The second moment of the distribution function per the unit of
mass is given by the expression $\frac {1}{\rho}P_{ij}(t,X) =
\left \langle v_iv_j, f_\mu^0(X,V) \right \rangle$. This
expression defines the tensor of pressure $P_{ij}(t,X)$ or, by the
other words, the tensor of viscose tensions (Using one or another
term depends on what kind of limiting transfer was targeted for
consideration). This tensor compensates the lowering of the
temperature due to collisions at the cost of shear friction. The
preceding equation takes form

\beq
 \frac{\partial (\rho (t,X)u_j(t,X))}{\partial t} + \sum_{i=1}^2
\frac {\partial}{\partial x_i} (\rho (t,X) P_{ij}(t,X)) + \rho
(t,X) {\cal F}_j(X) = 0.
 \label{11}
 \eeq

By multiplication of (\ref{7}) by $(W)^2$ and by integration over
$V$ one obtains the equation for the kinetic temperature (\ref{6})

\beq
\begin{array}{c}
 2\frac{\partial (\rho (t,X) \kappa_B T(t,X))}{\partial t} +
 \int (V-U(t,X))^2 \sum_{i=1}^2
\frac {\partial f^0}{\partial x_i} v_i dV + \\
 \int (V-U(t,X))^2 \sum_{i=1}^2 \frac {\partial f^0}{\partial v_i}
  {\cal F}_i(X) dV = - \zeta T(t,X),
\end{array}
 \label{12}
 \eeq
where

$$
 \zeta = - \frac {1}{2\rho \kappa_B T}\int dV (V- U(t,x))^2J[f]
$$
is the rate of cooling induced by collisions of particles. Let us
introduce the notation for the vector of the heat flow, that
corresponds to the second summand in (\ref{12}): $q(t,X) =\int
(V-U(t,X))^2f^0(t,X,V)VdV$. Collision integral equals zero. Let us
carry out of the sign of the integral factor ${\cal F}_i(X)$ in
the equation (\ref{12}) and integrate by parts. The  outside of
the integral term is zero because $f^0$ relative to $V$ is
Gaussian and the remaining integral equals zero since it gives the
mean value of the chaotic velocity (the chaotic velocity itself
but not its square). The equation (\ref{12}) takes the form

\beq
 2\frac{\partial (\rho (t,X) \kappa_B T(t,X))}{\partial t} +
\frac {\partial }{\partial x_i}[\rho (t,X)\kappa_B q_i(t,X)]  =
 -\zeta T(t,X).
 \label{13}
 \eeq

The equation (\ref{13}) gives something like equation of
continuity for the balance of the heat.

Introduce the set of subregions $\Omega_{\delta} \subset \Omega$
with boundaries $\Gamma_\delta = \{R=R_1+\delta\}\cap
\{R=R_2-\delta\}$.

The last summand in (\ref{11}) tends to infinity as $R \to R_1$
and as $R \to R_2$, since the function $\rho (t,X)$ is bounded
from below by a positive constant. Consequently, at the steady
state, when $U_j(t,X)$ does not depend on $t$, we have

\beq
 \sum_{i=1}^2 \frac {\partial}{\partial x_i}(\rho (X) P_{ij}(X))
 \to \infty
 \label{14}
 \eeq
 \noindent as $\delta \to 0$ for any indices $j$.

Integrate the relation (\ref{14}) over any subregion
$\Lambda_\delta \subset \Omega$ which contains the part
$\Gamma_\delta$ of the boundary of $\Omega_\delta$ for any fixed
value of the index $j$. In the stationary situation the flux of
the vector field $P_{ij}$ through $\Gamma_\delta$ directed by the
interior normal tends to $+\infty$ as $\delta \to 0$. This flux
define forces induced by the pressure $P_{ij}$ that acts on
particles of the Ring. The compensation of the pressure would
require infinite velocities that is physically impossible, to say
nothing of the observed data. The resultant forces at the vicinity
of the boundary of the Ring are directed inward the Ring.
Particles of $\Gamma$ must to fly out loosing the gravitational
connection with Saturn or to strive toward the interior of the
Ring. As a result, the density of particles on a vicinity of the
border will rapidly decrease. Consequently, as in the case when
collisions do not take into account, under the influence of forces
of the pressure, the Ring began contract to its middle part. The
state of the disc cannot remain stationary.

Hence the steady state surface density in two-dimensional model
$B$ cannot be smooth relative to the space variables.

$\Box$

\bigskip

\begin{corollary}

The theorem 2 remains valid for any sufficiently wide sub-ring
which is separated from other parts of the Ring by two divisions
deprived of gravitating particles.
\end{corollary}

Theorems 1 and 2 are theorems about non-existence of the expected
standard solution. This fact follow that to find smooth positive
function that gives steady state stable solution of models $A$ and
$B$ is impossible. By the other word these solutions must be
unexpected and nonstandard.

We arrive to the necessity to consider irregular distribution
functions.

\bigskip

\section{The theory of resonances and \\
the concept of fractal behavior}

\bigskip

The idea of gravitational resonance is due to Maxwell. He study
oscillations of particles of the Ring and distinguished the proper
oscillations of the system and the forced ones, evoked by external
effects. He did not correlate the external forces with concrete
planets or satellites of Saturn. Maxwell considered resonances as
phenomena of prime importance. Later on it was discovered one of
such resonances: it takes 22.6 hours to the external part of the
sub-ring B to make a complete turn, that exactly coincides with
the period of rotation of Mimas --- the Moon nearest to the Ring
of Saturn. Moreover, the farthest point of the external part of
this sub-ring (his apo-centrum) always consists the angle $\pi/4$
with the direction to Mimas. After this discovery, astronomers
tried to bound any sufficiently big gap with his own ``shepherd'',
i.e. with a satellite having the commensurable with the gap period
of rotation(to cause a resonance).

The term shepherd is a beautiful and bright poetical pattern. It
presuppose that each big satellite takes care of his own flock of
meteorites that constitute an isolated sub-ring. Admittedly, this
approach allows to predict positions of some small satellites of
Uranus by using irregularities of gaps in the Uranian Ring. The
discovery of the corresponding satellites of Uranus was the real
triumph of the concept of shepherds. But as regard to a real
physical sense of the term, it seems to us not fully adequate. It
is evident that each big satellite affects not solely on a one
sub-ring. All the set of satellites (the polyhedron of satellites)
affects on each sub-ring. All this structure as a whole is
responsible for all microstructure of the Ring.

Friedman and Gorkavyj \cite{GF} explored the Boltzmann equation
for the Ring. They took as a pattern the method of expanding of
solution in moments of velocity distribution of particles relative
to the stationary Maxwell distribution. This method was used by
Chapman end Enskog for investigation of plasma equations. In the
work \cite{GF} It was written out the kinetic Boltzmann equation
for three-dimensional model of the Ring, consisted of rigid
nonelastic particles executed the differential rotation in the
gravitational field of Saturn. For simplification of explorations,
it was considered a flat (two-dimensional) model which was an
analog to the model $B$ and the transport equations for it was
written.

Friedman and Gorkavyj certify, that solutions for these equations
are unknown. So, they consider the linearization of the equations
in the vicinity of Maxwell distribution. It is called by
dispersing equation. Further they consider various types of
oscillations of the disc density that meet this differential
equation. The distinct types of oscillations were obtained after
the neglect of one or another member of the equation with the
corresponded physical justifications.

\bigskip

The substantial interest represent oscillations of accretion type
bounded with the accretion of particles on the Ring from the
ambient space. This process supplements the diminishing of
particles due to falling some of them on Saturn. Indeed, as a
result of collision may appears the loss in the particle velocity
and such particles will try to fall on Saturn. But it should be
taken into account the following: if one accepts our conjecture of
the fractal structure of the Ring, than very thin and extremely
intensive jets of meteorites would prevent the falling of
particles. Many particles intended to fall on Saturn, being the
subject to impulses from jets, get trapped and tangled into the
Ring structure.

\bigskip

Conditions for instability for each of selected types of
oscillations were found in the work \cite{GF} . The necessity of
stratification of the Ring on the variety of thin ringlets was
explained by these instabilities. What kind of these conditions
appears fulfilled in the real situation remains unknown. It is not
so bad. It may be that such oscillations are realized in the Ring.
The real failure is the following: attempts of the strict
justification of its existence in the work \cite{GF} are
incorrect. The writing of variational equations as such presuppose
the differentiability of the equation (\ref{7}) relative to the
initial data. However, it was proved in the theorem 2 that the
equation (\ref{7}) can not have stationary, differentiable
solutions. Hence, all the written out oscillating modes were
obtained from internally inconsistent premises.

This failure is typical for many works. Using calculations
obtained from two-dimensional model and meeting obstacles, one
gets out of difficulty by appealing to three-dimensional real
situation. This position is contradictory. In contrast of this, we
explicitly declare that the two-dimensional model does not have
smooth solutions. The exit is to consider only three-dimensional
model, or to seek non-smooth (may be a fractal) solutions of the
two-dimensional model

Let us complicate the model by adding the impact of the Saturn
moons on particles of the Ring. Let $Z_i(t), \; i=1...l$ be
Cartesian coordinates of moons at a moment $t$ and $\nu_i$ being
its masses. The resultant gravitational force acting on a point $X
\in Q$ equals

$$
\Phi (t,X) = {\cal F}(X) + \sum_{i=1}^l\frac
{\nu_i}{|Z_i(t)-X|^2}.
$$

The Boltzmann equation takes the form

 \beq
 \begin{array}{c}
 \frac{\partial f_\mu(t,X,V)}{\partial t} +
 \sum_{i=1}^2\frac{\partial f_\mu(t,X,V)}{\partial x_i}v_i +
 \sum_{i=1}^2\frac{\partial f_\mu(t,X,V)}{\partial v_i}\Phi_i (t,X) = \\

  \frac {1}{\lambda}\int_0^\infty s(m)dm\int K(\sigma)
 [\tilde f_\mu \tilde f_m - f_\mu f_m]d\omega_\mu d\sigma.
 \end{array}
 \label{15}
 \eeq

\bigskip

{\it Our point of view is the following: The genuine reason for
stratification of the Ring consists in the fractal solution of the
Boltzmann equation connected with the piece-wise smooth
Hamiltonian system that is the limit-system for that of
(\ref{15}). The fractal solution is the optimal solution to the
affine in control extremal problem of minimization of the
functional of the action function. The part of the control play
gravitational forces from the Moons of Saturn $\frac
{\nu_i}{|Z_i(t)-X|^2}$. The sum of these forces takes its values
in the polyhedron with vertices in positions of the Moons
$Z_i(t)$.}

\bigskip

Let us explain. Long ago it has been established that the basis
for the lows of movement are variational principles. This concept
allows to introduce the theory of the Saturn Ring into the
framework of the optimal control theory and of the variational
problems.

Particles of the Ring moves, in response to the gravitational
force of Saturn, mainly by circles, creating the structure of the
Ring. Hence, in view of stationarity of the movement, it shows the
phase portrait of trajectories emanating practically from all
initial points of the Ring. By other words, we see the projection
on the configuration space of the optimal synthesis for the
minimization problem of the action function.

\bigskip

After Poincar\'e, the exploration of the structure of trajectories
of ordinary differential equations (ODE) with smooth right hand
size is based on the analysis of basic singularities of model
system: its singular points, limit cycles, attractors and so on.
After that, one proves that the phase portrait, being perturbed,
remains topologically equivalent to that of the initial model
system. Hamiltonian systems with discontinues right hand side
arising from Pontryagin Maximum Principle have specific structure
(so called tangential jumps). In this case the specific phenomenon
of general position takes place: the infinite number of switchings
on  finite intervals of trajectories (I Kupka \cite{K},
M.I.Zelikin and V.F.Borisov \cite{Z1}). In \cite{Z1} the theorem
on fiber bundle was proved. Roughly speaking, it claims:

Perturbations of Fuller's problem by members of higher order
relative to the scale group and by adding rectifiable auxiliary
variables do not change the topological structure of Fuller
problem synthesis (in particular, the infinite number of
switchings on finite intervals of trajectories). The neighborhood
of singular extremals of second order is fibred by manifolds
having perturbed Fuller's problem synthesis.

Hence, it was proved that the synthesis of the Fuller problem is a
model for Pontryagin-type Hamiltonian systems with one-dimensional
control. The generalizations  of this construction for problems
with one-dimensional control in infinite-dimensional spaces were
obtained in \cite{Z2}. In case of multi-dimensional control it was
found modes, where the optimal control performs in a finite time
the countable number of revolutions, and also the modes where the
optimal control accomplishs in a finite time the full circuit
along the everywhere dense winding of a torus \cite{Z3}.

This approach was extended to problems with the affine
multidimensional control when control variables take values in a
polyhedron. It was discovered a new general phenomenon for
Pontryagin-type Hamiltonian systems with discontinues right hand
side --- the stochastic dynamics --- the full circuit of
Cantor-like non-wandering points which is realized in a finite
time \cite{Z4}.

\bigskip

Naturally, one cannot apply this theorem directly to the behavior
of particles of the Saturn Ring.

At the first glance the theorem on piece-wise smooth Hamiltonians
is inapplicable to our case. Hamiltonian of the many body problem
is smooth. However, if one considers the action of several bodies
on a test particle, then in the process of moving arise situations
when the Hamiltonian appears as close as desired to a piece-wise
smooth one. Though linear combinations of gravitational forces
from Moons lies strictly inside of the polyhedron with vertices in
positions of Moons and the set of controls is an open polyhedron,
however, optimal trajectories for the problem with the open
polyhedron can be approximated by those of closed ones. But the
corresponding passage to the limit for the entire optimal
synthesis demands justifications.

Besides, the theorem of Zelikin-Lokutsievskiy-Hildebrand relates
to the ordinary differential equations. The Boltzmann equation is
integro-differential and infinite-dimensional. Nevertheless, the
more reason is to expect the fractal structure of its solutions in
projection of an infinite-dimensional picture on the
finite-dimensional configuration space.

\bigskip

The key idea is the following: On particles act external (relative
to the Ring and Saturn) forces -- that from numerous Moons of
Saturn. Although one can not dispose its positions, nevertheless,
they can be regarded as control, moreover, as the optimal control.
Indeed, the Moons together with Saturn and its Ring as a whole,
behave as if it were their aim to minimize the functional of
action. In this process the control being a linear combination of
forces applied at the finite number of points $Z_i(t)$
(corresponding to disposition of Moons) changes in a polyhedron.

\bigskip

\begin{conjecture}
Kinetic Boltzmann equation related to the system (\ref{15}) has a
fractal solution.
\end{conjecture}

To use an abstract mathematical theorem in a real physical
situation demand caution. The exact fractal solution may be
observed if the asteroids constituent the Ring were ideal points.
But as far as they have a finite size, although insignificantly
small in comparison with the magnitude of the Ring itself, one
observes only the approximative picture: the tremendous number
inserted one into another isolated ringlets consisted of particles
of finite size. This matter is too rough to realize theoretical
Cantor-like set of orbits. The situation call to mind that one
which took place while the Lorenz attractor was discovered: the
better were optical instruments the bigger number of ringlets one
could see. Particles of the Ring do its best (as far as it would
be possible for particles of centimetre-sized diameters) to
reproduce the exact fractal solution that minimizes the functional
of action. Different oscillating processes in the Ring are simply
vibrations around the faithful fractal solution.

\bigskip

The influence of the electromagnetic factors on the Ring structure
is expected explorations and explanations. This influence must be
considerable. The electric impulses, that were accepted by
interplanetary apparatus Voyager with the source being localized
in the region of the Ring, were 100000 times as large as the most
powerful lightnings in the Earth atmosphere. The radial strips, so
called ``spokes'', observed on the Ring, rotate as a distinct
units, retaining its straight form, independently of the
differential rotation of the whole ambient ``wheel'', with the
speed of rotation of the magnetic field of Saturn. The spacecraft
Voyager detected the flow of very thin charged particles (with the
size of approximately $10^{-8}$ m) having the speed close to 100
km/sek. It is natural to suppose that, due to the shear friction,
the differential rotation of particles must induce the ionization
of the mass of very thin dust and the accumulation of very big
charges in view of the vacuum isolation. As far as the Boltzmann
equation is applicable both to the kinetic effects and to the
behaviour of plasma, it seems to be an appropriate technique for
these explorations.

\bigskip

\section{Concluding remarks}

\bigskip

Let us remark, that all discovered Rings of planets (around
Jupiter, Saturn, Uranus, and  Neptune) consist of the set of
inserted one into another ringlets. Besides, the Solar System
itself has a Ring --- the Asteroid Belt situated between orbits of
Mars and Jupiter.

Our experience in exploration of the topological structure of
phase portraits and the proved fact about the fractal structure of
optimal syntheses in general situation, gives grounds to formulate
the following conjecture:

\begin{conjecture}

The Asteroid Belt of the Solar System approximates the exact
fractal solution of the corresponding variational problem and
consists of the very big number isolated, inserted one into
another ringlets. The part of shepherds play planets of Solar
System --- Mars, Jupiter, Saturn, Uranus...
\end{conjecture}

Surely, all the planet of the Solar System define corresponded
mobile polyhedron; its gravitational forces realize control of the
Asteroid Belt. These forces should be considered not as external
impacts that call resonances but as an intrinsic forces of the
full system minimizing the action function. That must be the cause
of fractal structure of the Asteroid Belt.

Gaps in the Asteroid Belt have been discovered. They are called
the Kirkwood gaps. Periods of its rotation until now were
associated with that of Jupiter. Admittedly, Mars has
significantly less mass than Jupiter but in the course of its
rotation Mars practically intersects the Asteroid Belt. In
addition the impact of other planet should be worth to take into
account.

We know that the big planets of the Solar System, and the Solar
System itself have Rings. So the right statement of the question
is not to ask: why planets have rings, but to ask: why the planets
of the Earth group (Mercury, Venus, Earth, and Mars) are devoid of
Rings. The standard preliminary conjecture is that they are too
near to Sun, where the density of the interplanetary dust is too
poor to complete its reduction due to accretion on planets. May be
rings do not sustain against the pressure of the Sun wind, or the
Sun light?

But in accordance with the our conception the cause can include in
the dimension of control. The theorem of
Zelikin-Lokutsievskiy-Hildebrand \cite{Z4}, \cite{Z5} demands that the dimension of control must be greater than 1. Mars has only two
satellites Fobos and Deimos, and the control is one-dimensional.
Earth has only one satellite. Venus and Mercury do not have
satellites. To the contrary, Jupiter, Saturn, Uranus, and Neptune
(all of them have Rings) have many satellites. Our point of view
is the following: satellites do not only stratify rings. The
creation and the maintenance of the existence of Rings are
connected with the presence of several big satellites. They
organize streams of meteorites, gathering its into a Ring, and
maintain the existence of the Ring by giving to it the fractal
structure.

We insist on the following proposition:

\bigskip

{\sl Existence of Rings and its stratification on ringlets is not
the exclusive but the regular phenomenon that is inherent in
sufficiently large planet systems}

\bigskip

Let us formulate the following program that seems to us very
important and very difficult as for the theory of kinetic
Boltzmann equation itself, so also for its numerous applications
in kinetic theory of gases, in cosmogony, in dynamics of plasma
and in other domains.

\bigskip

\begin{program}

 Investigation of fractal solutions to kinetic Boltzmann equation.
\end{program}
\bigskip

The first step to realize this program would be to explore the
above-mentioned simplest case and to prove of the conjecture 1.

We only have to say that on the way of this proof stands
considerable mathematical difficulties.  Questions of the
existence, of the smoothness, and of the extendability of
solutions to the Boltzmann equation evoke serious discussions
among professionals \cite{D},\cite{Ca} \cite{Po} \cite{U}. Since
the proof of the existence of the fractal solution to ordinary
differential equations appears very difficult (the exact and thorough proof takes about hundred pages), the analog of the
theorem for the integro-differential equation must be much more
difficult. But even this would not solve the problem because the
Boltzmann equation itself was derived under the assumption that
the density function is smooth.

Both Maxwell and Poincar\'e refused (not without reason) to make
final conclusions about the stability of Saturn Ring, although
namely this question was the subject of their searchings.
Moreover, Maxwell was led to conclude that the Ring will evolute
and eventually disintegrates in the foreseeable future. After
receiving from Struve observatory some not so high-reliable
information about changes in the size of the Saturn Ring, Maxwell
wrote:

"It will be worth while to investigate more carefully whether
Saturn's Rings are permanent or transitionary elements of the
Solar System, and whether in that part of the heavens we see
celestial immutability or terrestrial corruption and generation,
and the old order giving place to new before our own eyes".

\end{document}